\documentclass[12pt,draftcls,onecolumn]{IEEEtran}
\usepackage{amsmath,amssymb}
\usepackage{epstopdf}
\usepackage{latexsym}
\usepackage{psfrag}
\usepackage{epsfig}
\usepackage{graphicx,subfigure}
\usepackage{graphicx}
\usepackage{bm}

\newtheorem{theorem}{Theorem}
\newtheorem{proposition}{Proposition}
\newtheorem{corollary}{Corollary}
\newtheorem{definition}{Definition}
\newtheorem{lemma}{Lemma}
\newtheorem{remark}{Remark}
\newtheorem{example}{Example}

\newcommand{\R}{{\rm I\!R}}

\begin{document}

\title{Products of Generalized Stochastic Sarymsakov Matrices}

\author{Weiguo Xia, Ji Liu, Ming Cao, Karl H. Johansson, and Tamer Ba\c{s}ar%
\thanks{Weiguo Xia and Karl H. Johansson are with ACCESS Linnaeus Centre, School of Electrical Engineering, Royal Institute of Technology, Sweden.
 {\tt\small \{weiguo,kallej\}@kth.se}. Ji Liu and Tamer Ba\c{s}ar are with the Coordinated Science Laboratory, University of Illinois at Urbana-Champaign, USA. {\tt\small \{jiliu,basar1\}@illinois.edu}. Ming Cao is with the Faculty of Mathematics and Natural Sciences, ENTEG, University of Groningen, The Netherlands. {\tt\small m.cao@rug.nl}. }}

\maketitle

\begin{abstract}
In the set of
stochastic, indecomposable, aperiodic (SIA) matrices, the class of
stochastic Sarymsakov matrices is the \emph{largest} known subset (i) that
is closed under matrix multiplication and (ii) the infinitely long
left-product of the elements from a compact subset converges to a rank-one matrix.
In this paper, we show that a larger subset with these two properties can be derived  by
generalizing the standard definition for Sarymsakov matrices. The
generalization is achieved either by introducing an ``SIA index",
whose value is one for Sarymsakov matrices, and then looking at
those stochastic matrices with larger SIA indices, or by considering matrices that are not even SIA. Besides constructing a larger set, we give sufficient conditions for generalized
Sarymsakov matrices so that their products converge to rank-one
matrices. The new insight gained through studying generalized
Sarymsakov matrices and their products has led to a new understanding
of the existing results on consensus algorithms and will be helpful
for the design of network coordination algorithms.

\end{abstract}


\section{Introduction}

Over the last decade, there has been considerable interest in
consensus problems
\cite{Ts84,TsBeAt86,JaLiMo03,OlMu04,Mo05b,ReBe05,CaMoAn08a,ToNe14}
that are concerned with a group of agents trying to agree on a
specific value of some variable. Similar research interest arose
decades ago in statistics \cite{De74}.
While different aspects of consensus processes, such as convergence
rates \cite{OlTs09}, measurement delays \cite{CaMoAn08b}, stability
\cite{Mo05b,LiMoNeBa14}, and controllability \cite{EgMaCaCaBi12}, have been investigated, and many variants of
consensus problems, such as average consensus \cite{DiKaMoRaSc10},
asynchronous consensus \cite{CaMoAn08b}, quantized consensus
\cite{KaBaSr07}, and constrained consensus \cite{NeOzPa10}, have been
proposed,
some fundamental issues of discrete-time linear consensus processes
still remain open.

A discrete-time linear consensus process can typically be modeled by a linear recursion equation
of the form
\begin{equation}\label{eq:system}
x(k+1)=P(k)x(k),\ k\geq1,
\end{equation}
where $x(k)=[x_1(k),\ldots,x_n(k)]^T\in\R^n$ and each $P(k)$ is an
$n\times n$ stochastic matrix. It is well known that reaching a
consensus for any initial state in this model is equivalent to the product  $P(k)\cdots P(2)P(1)$ converging to a rank-one
matrix as $k$ goes to infinity. In this context, one fundamental
issue is as follows. Given a set of $n\times n$ stochastic matrices
$\mathcal{P}$, what are the conditions on $\mathcal{P}$ such that
for any infinite sequence of matrices $P(1),P(2),P(3),\ldots$ from
$\mathcal{P}$, the sequence of left-products
$P(1),P(2)P(1),P(3)P(2)P(1),\ldots$ converges to a rank-one matrix?
We will call $\mathcal{P}$ satisfying this property a consensus set
(the formal definition will be given in the next section). The
existing literature on characterizing a consensus set can be traced
back to at least the work of Wolfowitz \cite{Wo63} in which
stochastic, indecomposable, aperiodic (SIA) matrices are introduced.
Recently, it has been shown in \cite{BlOl14} that deciding whether
$\mathcal{P}$ is a consensus set is NP-hard; a combinatorial necessary
and sufficient condition for deciding a consensus set has also been
provided there. Even in the light of these classical and recent
findings, the following fundamental question remains: What is the
largest subset of the class of $n\times n$ stochastic matrices whose compact
subsets are all consensus sets? In \cite{LiMoAnYu11}, this question is
answered under the assumption that each stochastic matrix has positive diagonal entries. For general stochastic matrices,
the question has remained open. This paper aims at
dealing with this challenging question by checking some well-known classes of SIA matrices.

In the literature, the set of stochastic
Sarymsakov matrices, first introduced by Sarymsakov \cite{Sa61}, is
the \emph{largest} known subset of the class of stochastic matrices whose compact
subsets are all consensus sets; in particular, the set is closed
under matrix multiplication and the left-product of the elements from its compact subset
converges to a rank-one matrix \cite{Ha02}. In
this paper, we construct a larger set of stochastic matrices whose compact subsets are all
consensus sets. We adopt the natural idea which is to generalize the
definition of stochastic Sarymsakov matrices so that the original
set of stochastic Sarymsakov matrices are contained.

In the paper, we  introduce two ways to generalize the definition
and thus study two classes of generalized stochastic Sarymsakov
matrices. The first class makes use of the concept of the SIA index
(the formal definition will be given in the next section). It is shown that the set of $n\times n$ stochastic matrices
with SIA index no larger than $k$ is closed under matrix
multiplication only when $k=1$, which turns out to be the stochastic
Sarymsakov class; this result reveals why exploring a consensus set
larger than the set of stochastic Sarymsakov matrices is a challenging
problem. A set that consists of all stochastic
Sarymsakov matrices plus one specific SIA matrix and thus is
slightly larger than the stochastic Sarymsakov class is constructed,
and we show that it is closed under matrix multiplication. For the other
class of generalized Sarymsakov matrices, which contains matrices
that are not SIA, sufficient conditions are provided for the convergence of the product
of an infinite matrix sequence  from this class to a
rank-one matrix. A special case in which all the generalized
Sarymsakov matrices are doubly stochastic is also discussed.

The rest of the paper is organized as follows. Preliminaries are introduced in Section \ref{se:pre}. Section \ref{se:index} introduces the SIA index and discusses the properties of the set of stochastic matrices with SIA index no larger than $k,\ k\geq1$. In Section \ref{se:W}, sufficient conditions are provided for the convergence of the product of an infinite sequence of matrices from a class of generalized stochastic Sarymsakov matrices, and the results are applied to the class of doubly stochastic matrices.  Section \ref{se:conclusion} concludes the paper.


\section{Preliminaries}\label{se:pre}

We first introduce some basic definitions. Let $n$ be  a positive
integer. A square matrix $P=\{p_{ij}\}_{n\times n}$ is said to be \emph{stochastic} if $p_{ij}\geq0$ for all
$i,\ j \in \{1,\ldots,n \}=\mathcal N$, and $\sum_{j=1}^np_{ij}=1$ for all
$i=1,\ldots,n$. Consider a stochastic matrix $P$. For a set
$\mathcal{A}\subseteq\mathcal N$, the
set of \emph{one-stage consequent indices} \cite{Se79} of
$\mathcal{A}$ is defined by
$$F_P(\mathcal{A})=\{j:\ p_{ij}>0\ \text{for
some}\ i\in \mathcal{A}\}$$
and we call $F_P$ the consequent function of $P$. For a singleton $\{i\}$, we write $F_P(i)$ instead of $F_P(\{i\})$ for simplicity. A matrix $P$ is indecomposable and aperiodic, and thus called an \emph{SIA} matrix,
if $\lim_{m\rightarrow\infty}P^m=\bm{1}c^T,$ where $\bm{1}$
is the $n$-dimensional all-one column vector, and
$c=[c_1,\ldots,c_n]^T$ is some column vector satisfying $c_i\geq0$
and $\sum_{i=1}^nc_i=1$. $P$ is said to belong to the \emph{Sarymsakov
class} or equivalently $P$ is a {\em Sarymsakov} matrix if for any two disjoint nonempty sets
$\mathcal{A},\ \tilde{\mathcal{A}}\subseteq\mathcal N$, either
\begin{equation}\label{eq:Sarym1}
F_P(\mathcal{A})\cap F_P(\tilde{\mathcal{A}})\neq\emptyset,
\end{equation} or
\begin{equation}\label{eq:Sarym2}F_P(\mathcal{A})\cap F_P(\tilde{\mathcal{A}})=\emptyset \text{ and }
|F_P(\mathcal{A})\cup
F_P(\tilde{\mathcal{A}})|>|\mathcal{A}\cup\tilde{\mathcal{A}}|,\end{equation}
where $|\mathcal{A}|$ denotes the cardinality of  $\mathcal{A}$. We
say that $P$ is a \emph{scrambling} matrix if for any pair of distinct
indices $i,j\in\mathcal N$, $F_P(i)\cap F_P(j)\neq \emptyset$, which is equivalent to requiring that there always exists an index $k\in\mathcal N$ such
that both $p_{ik}$ and $p_{jk}$ are positive.

From the definitions, it should be obvious that a scrambling matrix belongs to the Sarymsakov class. It has been proved in \cite{Se79}
that any product of $n-1$ matrices from the Sarymsakov class is scrambling. Since a stochastic scrambling matrix is SIA \cite{Se06}, any stochastic Sarymsakov matrix must be an SIA matrix.

\begin{definition}({\it Consensus set})\label{def:consensus}
Let $\mathcal{P}$ be a set of $n\times n$ stochastic matrices. $\mathcal P$ is a consensus set if for each sequence of matrices $P(1),\ P(2),\ P(3),\ldots$ from $\mathcal{P}$, $P(k)\cdots P(1)$ converges to a rank-one matrix $\bm{1}c^T$ as $k\rightarrow\infty$, where $c_i\geq0$
and $\sum_{i=1}^nc_i=1$.
\end{definition}

Deciding whether a set is  a consensus set or not is critical in establishing the convergence of the state of system (\ref{eq:system}) to a common value. Necessary and sufficient conditions for $\mathcal{P}$ being a consensus set have been established \cite{Wo63,AnTi77,Se06,BlOl14,XiCa14}.

\begin{theorem}\cite{XiCa14}\label{thm:products1}
Let $\mathcal{P}$ be a compact set of $n\times n$ stochastic  matrices. The following conditions are equivalent:
\begin{enumerate}
\item $\mathcal P$ is a consensus set.
\item For each integer $k\geq1$ and any $P(i)\in\mathcal{P},\ 1\leq i\leq k$, the matrix $P(k)\cdots P(1)$ is SIA.
\item There is an integer $\nu\geq1$ such that for each $k\geq \nu$ and any $P(i)\in\mathcal{P},\ 1\leq i\leq k$, the matrix $P(k)\cdots P(1)$ is scrambling.
\item There is an integer $\mu\geq1$ such that for each $k\geq \mu$ and any $P(i)\in\mathcal{P},\ 1\leq i\leq k$, the matrix $P(k)\cdots P(1)$ has a column with only positive elements.
\item There is an integer $\alpha\geq1$ such that for each $k\geq \alpha$ and any $P(i)\in\mathcal{P},\ 1\leq i\leq k$, the matrix $P(k)\cdots P(1)$ belongs to the Sarymsakov class.
\end{enumerate}
\end{theorem}

For a compact set $\mathcal P$ to be a consensus set, it is necessary that every matrix in $\mathcal P$ is SIA in view of item (2) in Theorem \ref{thm:products1}. If a set of SIA matrices is closed under matrix multiplication, then one can easily conclude from item (2) that its compact subsets are all consensus sets. However, it is well known that the product of two SIA matrices may not be SIA. The stochastic Sarymsakov class is the largest known set of stochastic matrices, which is closed under matrix multiplication. Whether there exists a larger class of SIA matrices, which contains the Sarymsakov class as a proper subset and is closed under matrix multiplication, remains unknown. We will explore this by taking a closer look at the definition of the Sarymsakov class and study the properties of classes of generalized Sarymsakov matrices that contain the Sarymsakov class
as a subset.


\section{SIA index}\label{se:index}

The key notion in the definition of the Sarymsakov class is the set of one-stage consequent indices. We next introduce the notion of the set of {\it $k$-stage consequent indices} and utilize this to define a larger matrix set, which contains the Sarymsakov class.

For a stochastic matrix $P$ and a set $\mathcal{A}\subseteq\mathcal N$, let $F_P^k(\mathcal{A})$ be
the set of {\it $k$-stage consequent indices} of $\mathcal{A}$,
which is defined by
$$F_P^1(\mathcal{A})=F_P(\mathcal{A})\ \mbox{and } F_P^k(\mathcal{A})=F_P(F_P^{k-1}(\mathcal{A})),\ k\geq2.$$

\begin{definition}A stochastic matrix $P$ is said to belong to the class $\mathcal{S}$ if
for any two disjoint nonempty subsets $\mathcal{A},\
\tilde{\mathcal{A}}\subseteq\mathcal N$, there exists an integer
$k\geq1$ such that either
\begin{equation}\label{eq:SIA1}F_P^k(\mathcal{A})\cap
F_P^k(\tilde{\mathcal{A}})\neq\emptyset,\end{equation}
or
\begin{equation}\label{eq:SIA2}F_P^k(\mathcal{A})\cap
F_P^k(\tilde{\mathcal{A}})=\emptyset\ \text{and}\
|F_P^k(\mathcal{A})\cup
F_P^k(\tilde{\mathcal{A}})|>|\mathcal{A}\cup\tilde{\mathcal{A}}|.\end{equation}
\end{definition}

It is easy to see that the Sarymsakov class is a subset of $\mathcal{S}$ since $k$ is 1 in the definition of  the Sarymsakov class. An important property of the consequent function $F_P$ given below will be useful.

\begin{lemma} \cite{Ha02} \label{lm:consequent} Let $P$ and $Q$  be $n\times n$ nonnegative matrices. Then,
$F_{PQ}(\mathcal{A})=F_{Q}(F_P(\mathcal{A}))$ for all subsets
$\mathcal{A} \subseteq\mathcal N$.
\end{lemma}

A direct consequence of Lemma \ref{lm:consequent} is that $F_{P^k}(\mathcal{A})=F_{P}^k(\mathcal{A})$ for any stochastic matrix $P$, any integer $k\geq1$ and any subset $\mathcal A\subseteq\mathcal N$.

The following theorem establishes the relationship between the matrices in $\mathcal S$ and SIA matrices.

\begin{theorem} \cite{XiCa14}\label{th:SIA1} A stochastic matrix $P$ is in $\mathcal{S}$ if and
only if $P$ is SIA.
\end{theorem}

From Theorem 4.4 in \cite{Pa71}, we know the following result.
\begin{theorem} \cite{Pa71}\label{th:SIA2} A stochastic matrix $P$ is SIA if and only if for every pair of indices $i$ and $j$, there exists an integer $k$, $k\leq n(n-1)/2$, such that $F^k_P(i)\cap F^k_P(j)\neq\emptyset$.
\end{theorem}

Theorem \ref{th:SIA2} implies that the index $k$ in (\ref{eq:SIA1}) and (\ref{eq:SIA2}) can be bounded by some integer.

\begin{lemma}\label{lm:SIAorder2}
A stochastic matrix $P$ is SIA if and only if for any pair of disjoint nonempty sets $\mathcal{A},\ \tilde{\mathcal{A}}\subseteq\mathcal N$, there exists an index $k$, $k\leq n(n-1)/2$, such that $F^k_P(\mathcal{A})\cap F^k_P(\tilde{\mathcal{A}})\neq\emptyset$.
\end{lemma}

\begin{example}
Let
\begin{equation}
P=\begin{bmatrix}
      \frac{1}{3} &\frac{1}{3}  & \frac{1}{3}  \\
      1 & 0 & 0\\
      0 & 1 & 0 \\
    \end{bmatrix}.
    \end{equation}
$P$ is a stochastic matrix. Consider two disjoint nonempty sets $\mathcal A=\{2\},\tilde{\mathcal{A}}=\{3\}$. One knows that $F_P(\mathcal{A})=\{1\}$ and $F_P(\tilde{\mathcal{A}})=\{2\}$, implying that $F_P(\mathcal{A})\cap F_P(\tilde{\mathcal{A}})\neq\emptyset$ and $|F_P(\mathcal{A})\cup F_P(\tilde{\mathcal{A}})|=|\mathcal A\cup\tilde{\mathcal{A}}|$. Therefore, $P$ is not a Sarymsakov matrix. However, the fact that $F^2_P(\mathcal{A})=\{1,2,3\}$ and $F^2_P(\tilde{\mathcal{A}})=\{2\}$ shows that $F^2_P(\mathcal{A})\cap F^2_P(\tilde{\mathcal{A}})\neq\emptyset$. This means (\ref{eq:SIA1}) holds for $k=2.$

For every other pair of disjoint nonempty sets $\mathcal A,\tilde{\mathcal{A}}\subseteq\mathcal N$, it can be verified that $F_P(\mathcal{A})\cap F_P(\tilde{\mathcal{A}})\neq\emptyset$. One has that though $P$ is not a Sarymsakov matrix, $P$ is an SIA matrix from  Lemma \ref{lm:SIAorder2}. \hfill $\Box$
\end{example}

From the above example and Lemma \ref{lm:SIAorder2}, one knows that the class of SIA matrices may contain a large number of matrices that do not belong to the Sarymsakov class. Starting from the Sarymsakov class, where $k=1$ in (\ref{eq:SIA1}) and (\ref{eq:SIA2}), we relax the constraint on the value of the index $k$ in (\ref{eq:SIA1}) and (\ref{eq:SIA2}), i.e., allowing for $k\leq2$, $k\leq 3,\dots,$ and obtain a larger set containing the Sarymsakov class. We formalize the idea below and study whether the derived set is closed under matrix multiplication or not.

Consider a fixed integer $n$. Denote all the unordered pairs of disjoint nonempty sets of $\mathcal N$ as $(\mathcal{A}_1,\ \tilde{\mathcal{A}}_1),\ldots,(\mathcal{A}_m,\ \tilde{\mathcal{A}}_m)$, where $m$ is a finite number.

\begin{definition}\label{def:index}
Let $P\in\R^{n\times n}$ be an SIA matrix. For each pair of disjoint nonempty sets $\mathcal{A}_i,\ \tilde{\mathcal{A}}_i\subseteq\mathcal N,\ i=1,\ldots,m$, let $s_i$ be the smallest integer such that  either (\ref{eq:SIA1}) or (\ref{eq:SIA2}) holds. The {\em SIA index} $s$ of $P$ is $s=\max\{s_1,s_2,\ldots,s_m\}$.
\end{definition}

From Lemma \ref{lm:SIAorder2}, we know that for an SIA matrix $P$, its SIA index $s$ is upper bounded by $n(n-1)/2$. Assume that the largest value of the SIA indices of all the $n\times n$ SIA matrices is $l$, which depends on the order $n$. We define several subsets of the class of SIA matrices. For $1\leq k\leq l,$  let
\begin{align}\label{eq:V}
\mathcal{V}_k&=\{P\in\R^{n\times n}| P \text{ is SIA and its SIA index is } k\}
\end{align}
and
\begin{align}\label{eq:SIAsub}
\mathcal{S}_k&=\cup_{r=1}^k \mathcal V_r.
\end{align}

Obviously $\mathcal{S}_1\subset\mathcal{S}_2\subset\cdots\subset\mathcal{S}_l$ and $\mathcal{S}_1=\mathcal V_1$ is the class of stochastic Sarymsakov matrices. One can easily check that when $n=2$, all SIA matrices are scrambling matrices and hence belong to the Sarymsakov class. When $n\geq3,$ the set $\mathcal V_{n-1}$ is  nonempty. To see this, consider an $n\times n$ stochastic matrix
$$P=\begin{bmatrix}
      \frac{1}{n} & \frac{1}{n}  & \cdots & \frac{1}{n}  & \frac{1}{n}  \\
      1 & 0 & \cdots & 0 & 0 \\
      0 & 1 & \cdots & 0 & 0 \\
      \vdots & \vdots & \ddots & \vdots & \vdots \\
      0 & 0 & \cdots & 1 & 0 \\
    \end{bmatrix}.$$
For an index $i\in\mathcal N,\ i\neq n$, it is easy to check that $F^{n-1}_P(i)=\mathcal N$. Hence, for any two nonempty disjoint sets $\mathcal A,\tilde{\mathcal A}\in\mathcal N$, it must be true that $F_P^{n-1}(\mathcal{A})\cap F_P^{n-1}(\tilde{\mathcal{A}})\neq\emptyset$, proving that $P$ is an SIA matrix. Consider the specific pair of sets $\mathcal A=\{n\},\tilde{\mathcal A}=\{n-1\}$. One has that $F_P^{n-2}(n)=\{2\}$, $F_P^{n-2}(n-1)=\{1\}$, and $F_P^{n-1}(n)\cap F_P^{n-1}(n-1)\neq\emptyset$, implying that $P\in\mathcal V_{n-1}$. From this example, we know that a lower bound for $l$ is $n-1$.

In the next three subsections, we first discuss the properties of $\mathcal S_i,\ i=1,\ldots,l$, then construct a set,
closed under matrix multiplication, consisting of a specific  SIA matrix and all stochastic Sarymsakov matrices, and finally discuss the class of pattern-symmetric matrices.

\subsection{Properties of $\mathcal S_i$}\label{se:Si}

The following novel result reveals the properties of the sets $\mathcal S_i,\ 1\leq i\leq l$.

\begin{theorem}\label{thm:Sarymlargest}
Suppose that $n\geq3$. Among the sets $\mathcal{S}_1,\ \mathcal{S}_2,\ldots,\ \mathcal{S}_l$, the set $\mathcal{S}_1$ is the only set that is closed under matrix multiplication.
\end{theorem}

Note that a compact subset $\mathcal P$ of $\mathcal S_1$ is a consensus set. However, if  $\mathcal P$ is a compact set consisting of matrices in $\mathcal V_i,\ i\geq2$, Theorem \ref{thm:Sarymlargest} shows that there is no guarantee that $\mathcal P$ is still a consensus set.

The proof of Theorem \ref{thm:Sarymlargest} relies on the following key lemma, based on which the conclusion of Theorem \ref{thm:Sarymlargest} immediately follows. Before stating the lemma, we define a matrix $Q$ in terms of a matrix $P\in\mathcal{V}_i,\ i\geq2$.

For a given matrix $P\in\mathcal{V}_i,\ i\geq2,$ from the definition of the Sarymsakov class, one has that there exist two disjoint nonempty sets $\mathcal{A},\ \tilde{\mathcal{A}}\subseteq \mathcal N$ such that $F_P(\mathcal{A})\cap F_P(\tilde{\mathcal{A}})=\emptyset$ and
\begin{equation}\label{eq:keylm1}
|F_P(\mathcal{A})\cup
F_P(\tilde{\mathcal{A}})|\leq|\mathcal{A}\cup\tilde{\mathcal{A}}|.\end{equation}
Define a matrix $Q=(q_{ij})_{n\times n}$ as follows
\begin{equation}\label{eq:keylm5}
q_{ij}=\begin{cases}
\frac{1}{|\mathcal{A}|}, & i\in F_{P}(\mathcal{A}),\ j\in\mathcal{A}, \\
0, & i\in F_{P}(\mathcal{A}),\ j\in\bar{\mathcal{A}}, \\
\frac{1}{|\tilde{\mathcal{A}}|}, & i\in F_{P}(\tilde{\mathcal{A}}),\ j\in\tilde{\mathcal{A}}, \\
0, & i\in F_{P}(\tilde{\mathcal{A}}),\ j\in\bar{\tilde{\mathcal{A}}}, \\
\frac{1}{n}, & \text{otherwise}, \\
\end{cases}
\end{equation}
where $\bar{\mathcal A}$ denotes the complement of $\mathcal A$ with respect to $\mathcal N$.

\begin{lemma}\label{lm:keylm}
Suppose that $n\geq 3$. For any $i=2,\ldots,l$, given a stochastic matrix $P\in\mathcal{V}_i$, then the matrix $Q$ given in (\ref{eq:keylm5}) belongs to the set $\mathcal{S}_2$ and $PQ,\ QP$ are not SIA. In addition, $Q\in\mathcal V_1$ if (\ref{eq:keylm1}) holds with the equality sign; $Q\in\mathcal V_2$ if the inequality (\ref{eq:keylm1}) is strict.
\end{lemma}

{\it Proof:} $Q$ is a stochastic matrix as each element $q_{ij},\ i,j=1,\ldots,n$ is nonnegative and each row  of $Q$ sums up to 1. Note that for an index $i\in\mathcal N$, the set of its one-stage consequent indices $F_Q(i)$ is either $\mathcal{A}$ or $\tilde{\mathcal{A}}$ or $\mathcal N$. We next show that $Q$ belongs to the set $\mathcal S_2$.

Consider two arbitrary disjoint nonempty sets $\mathcal{C},\ \tilde{\mathcal{C}}\subseteq\mathcal N$. One of the following statements must hold:
\begin{itemize}
\item[(a)] $\mathcal{C}\cup\tilde{\mathcal{C}}$ contains some element in $\overline{F_{P}(\mathcal{A})\cup F_{P}(\tilde{\mathcal{A}})}$;
\item[(b)] $\mathcal{C}\cup\tilde{\mathcal{C}}\subseteq F_{P}(\mathcal{A})\cup F_{P}(\tilde{\mathcal{A}})$, $\mathcal{C}\cap F_{P}(\mathcal{A})\neq\emptyset$, and $\tilde{\mathcal{C}}\cap F_{P}(\mathcal{A})\neq\emptyset$;
\item[(c)] $\mathcal{C}\cup\tilde{\mathcal{C}}\subseteq F_{P}(\mathcal{A})\cup F_{P}(\tilde{\mathcal{A}})$, $\mathcal{C}\cap F_{P}(\tilde{\mathcal{A}})\neq\emptyset$, and $\tilde{\mathcal{C}}\cap F_{P}(\tilde{\mathcal{A}})\neq\emptyset$;
\item[(d)] $\mathcal{C}\subseteq F_{P}(\mathcal{A}) \mbox{ and } \tilde{\mathcal{C}}\subseteq F_{P}(\tilde{\mathcal{A}});$
\item[(e)] $\mathcal{C}\subseteq F_{P}(\tilde{\mathcal{A}}) \mbox{ and }\tilde{\mathcal{C}}\subseteq F_{P}(\mathcal{A}).$
\end{itemize}

{\it Case (a).} From the definition of the matrix $Q$ in (\ref{eq:keylm5}), one has that $F_{Q}(\mathcal{C})$ or $F_{Q}(\tilde{\mathcal{C}})$ is the set $\mathcal N$, which implies that $F_{Q}(\mathcal{C})\cap F_{Q}(\tilde{\mathcal{C}})\neq\emptyset$.

{\it Case (b).} It is easy to see that $\mathcal{A}$ is a subset of both $F_{Q}(\mathcal{C})$ and $F_{Q}(\tilde{\mathcal{C}})$. Hence $F_{Q}(\mathcal{C})\cap F_{Q}(\tilde{\mathcal{C}})\neq\emptyset$.

{\it Case (c).} Similar to case (b), $\tilde{\mathcal{A}}$ is a subset of both $F_{Q}(\mathcal{C})$ and $F_{Q}(\tilde{\mathcal{C}})$. Hence $F_{Q}(\mathcal{C})\cap F_{Q}(\tilde{\mathcal{C}})\neq\emptyset$.

{\it Case (d).} From the definition of $Q$, one has
\begin{equation}\label{eq:keylm8}
F_{Q}(\mathcal{C})=\mathcal{A},\ F_{Q}(\tilde{\mathcal{C}})=\tilde{\mathcal{A}}.
\end{equation}
Following (\ref{eq:keylm1}),
$$|F_{Q}(\mathcal{C})\cup F_{Q}(\tilde{\mathcal{C}})|=|\mathcal{A}\cup\tilde{\mathcal{A}}|\geq|F_{P}(\mathcal{A})\cup F_{P}(\tilde{\mathcal{A}})|\geq|\mathcal{C}\cup\tilde{\mathcal{C}}|.$$
If $|F_{P}(\mathcal{A})\cup F_{P}(\tilde{\mathcal{A}})|>|\mathcal{C}\cup\tilde{\mathcal{C}}|$, then  $|F_{Q}(\mathcal{C})\cup F_{Q}(\tilde{\mathcal{C}})|>|\mathcal{C}\cup\tilde{\mathcal{C}}|$.

If $|F_{P}(\mathcal{A})\cup F_{P}(\tilde{\mathcal{A}})|=|\mathcal{C}\cup\tilde{\mathcal{C}}|$, we consider two cases
\begin{itemize}
\item[(d1)] $|\mathcal{A}\cup\tilde{\mathcal{A}}|>|F_{P}(\mathcal{A})\cup F_{P}(\tilde{\mathcal{A}})|$;
\item[(d2)] $|\mathcal{A}\cup\tilde{\mathcal{A}}|=|F_{P}(\mathcal{A})\cup F_{P}(\tilde{\mathcal{A}})|$.
\end{itemize}

{\it Case (d1).} We immediately conclude that $|F_{Q}(\mathcal{C})\cup F_{Q}(\tilde{\mathcal{C}})|>|\mathcal{C}\cup\tilde{\mathcal{C}}|$.

{\it Case (d2).} From
$$|\mathcal{A}\cup\tilde{\mathcal{A}}|=|F_{P}(\mathcal{A})\cup F_{P}(\tilde{\mathcal{A}})|=|\mathcal{C}\cup\tilde{\mathcal{C}}|,$$ one obtains  $\mathcal{C}=F_{P}(\mathcal{A})$ and $\tilde{\mathcal{C}}=F_{P}(\tilde{\mathcal{A}})$. We further look at  the sets of two-stage consequent indices of $\mathcal{C}$ and $\tilde{\mathcal{C}}$, and obtain from (\ref{eq:keylm8}) that
$$F^2_{Q}(\mathcal{C})=F_Q(\mathcal{A}),\ F^2_{Q}(\tilde{\mathcal{C}})=F_Q(\tilde{\mathcal{A}}).$$
We will show that $F_Q(\mathcal{A})\cap F_Q(\tilde{\mathcal{A}})\neq\emptyset$ from which we see that the smallest integer $k$ is 2 such that (\ref{eq:SIA1}) holds for this pair $\mathcal{C}$ and $\tilde{\mathcal{C}}$ and the matrix $Q$. Suppose on the contrary that $F_Q(\mathcal{A})\cap F_Q(\tilde{\mathcal{A}})=\emptyset$ is true. Since for any $i\in\mathcal N$, $F_Q(i)$ is either $\mathcal{A}$ or $\tilde{\mathcal{A}}$ or $\mathcal N$, the fact that $F_Q(\mathcal{A})\cap F_Q(\tilde{\mathcal{A}})=\emptyset$ implies that either
\begin{equation}\label{eq:keylm6}
F_Q(\mathcal{A})=\mathcal{A},\ F_Q(\tilde{\mathcal{A}})=\tilde{\mathcal{A}},
\end{equation}
or
\begin{equation}\label{eq:keylm7}
F_Q(\mathcal{A})=\tilde{\mathcal{A}},\ F_Q(\tilde{\mathcal{A}})=\mathcal{A}.
\end{equation}

If (\ref{eq:keylm6}) holds, then it is inferred from the structure of the matrix $Q$ that $\mathcal{A}\subseteq F_P(\mathcal{A})$ and $\tilde{\mathcal{A}}\subseteq F_P(\tilde{\mathcal{A}})$. Combining with the fact that $|F_P(\mathcal{A})\cup
F_P(\tilde{\mathcal{A}})|=|\mathcal{A}\cup\tilde{\mathcal{A}}|$, it must be true that $ F_P(\mathcal{A})=\mathcal{A}$ and $F_P(\tilde{\mathcal{A}})=\tilde{\mathcal{A}}$. It then follows that $F_P^k(\mathcal{A})=\mathcal{A}$ and $F_P^k(\tilde{\mathcal{A}})=\tilde{\mathcal{A}}$ for any positive integer $k$, showing that $P$ is not an SIA matrix in view of Lemma \ref{lm:SIAorder2}. We conclude that $F^2_{Q}(\mathcal{C})\cap F^2_{Q}(\tilde{\mathcal{C}})\neq\emptyset$.

If (\ref{eq:keylm7}) holds, then from the structure of the matrix $Q$ one has that $\mathcal{A}\subseteq F_P(\tilde{\mathcal{A}})$ and $\tilde{\mathcal{A}}\subseteq F_P(\mathcal{A})$. Similarly one obtains that $F_P(\mathcal{A})=\tilde{\mathcal{A}}$ and $F_P(\tilde{\mathcal{A}})=\mathcal{A}$. Thus, $F_P^k(\mathcal{A})\cap F_P^k(\tilde{\mathcal{A}})=\emptyset$  for any positive integer $k$, implying that $P$ is not an SIA matrix based on Lemma \ref{lm:SIAorder2}. We conclude that $F^2_{Q}(\mathcal{C})\cap F^2_{Q}(\tilde{\mathcal{C}})\neq\emptyset$.

{\it Case (e).} The discussion is similar to that in case (d).

Therefore, summarizing the discussions in all cases, we have shown that $Q\in\mathcal V_1$ if (\ref{eq:keylm1}) holds with the equality sign; $Q\in\mathcal V_2$ if the inequality (\ref{eq:keylm1}) is strict.

We next look at the matrix product $PQ$.  Consider the pair of sets $\mathcal{A}$ and $\tilde{\mathcal{A}}$. One has
\begin{equation}
F_{PQ}(\mathcal{A})=F_Q(F_P(\mathcal{A}))=\mathcal A,\ F_{PQ}(\tilde{\mathcal{A}})=F_{Q}(F_P(\tilde{\mathcal{A}}))=\tilde{\mathcal{A}}.
\end{equation}
Thus, for any positive integer $k$, $F^k_{PQ}(\mathcal{A})=\mathcal A$ and $F^k_{PQ}(\tilde{\mathcal{A}})=\tilde{\mathcal{A}}.$ $PQ$ is not an SIA matrix. Similarly, one can see that
\begin{equation}
F_{QP}(F_P({\mathcal{A}}))=F_P({\mathcal A}),\ F_{QP}(F_P(\tilde{\mathcal{A}}))=F_P(\tilde{\mathcal{A}}),
\end{equation}
implying that $QP$ is not an SIA matrix.
 \hfill $\Box$

\begin{remark}
Note that whether a stochastic matrix is SIA or not only depends on the positions of its nonzero elements but not their magnitudes. One can derive other matrices based on $Q$ in (\ref{eq:keylm5}) such that $PQ$ is not SIA by varying the magnitudes of the positive elements of $Q$ as long as each row sum equals to 1 and the positive elements are kept positive.
\end{remark}

There has been research work on defining an other class of matrices that is a subset of the SIA matrices and larger than the stochastic scrambling matrices. We establish its relationship with the stochastic Sarymsakov class in view of Lemma \ref{lm:keylm}.

\begin{definition}\cite{Se06}
$P\in \mathcal G,$ if $P$ is SIA and for any SIA matrix $Q$, $QP$ is SIA.
\end{definition}

\begin{proposition}
For $n\geq3,$ $\mathcal G$ is a proper subset of the class of stochastic Sarymsakov matrices $\mathcal S_1$.
\end{proposition}

{\it Proof.} Obviously $\mathcal G$ is a subset of $\mathcal S$. For any $P\in\mathcal V_i, \ i\geq2$, $P$ is not an element of $\mathcal G$ since there exists an SIA matrix $Q$ such that $QP$ is not SIA from Lemma \ref{lm:keylm}. Hence $\mathcal G$ is a subset of $\mathcal S_1$.

For $n\geq3$, let
$$P=\begin{bmatrix}
           \frac{1}{2} & \frac{1}{2} & 0 & 0 & \cdots & 0 \\
      0 & 0& 1 & 0 & \cdots & 0 \\
 \frac{1}{n} & \frac{1}{n} & \frac{1}{n} & \frac{1}{n} & \cdots   & \frac{1}{n} \\
\frac{1}{n} & \frac{1}{n} & \frac{1}{n} & \frac{1}{n} & \cdots &  \frac{1}{n} \\
      \vdots & \vdots & \vdots & \vdots & \ddots& \vdots \\
      \frac{1}{n} & \frac{1}{n} & \frac{1}{n}  & \frac{1}{n} & \cdots & \frac{1}{n} \\
    \end{bmatrix}.$$
One can verify that $P\in\mathcal S_1$. We show that $P\notin\mathcal G$. Consider the following matrix
$$Q=\begin{bmatrix}
           1 & 0 & 0 & 0 & \cdots & 0 \\
      1 & 0& 0 & 0 & \cdots & 0 \\
      0 & 1 & 0& 0 & \cdots & 0  \\
\frac{1}{n} & \frac{1}{n} & \frac{1}{n} & \frac{1}{n} & \cdots &  \frac{1}{n} \\
      \vdots & \vdots & \vdots & \vdots & \ddots& \vdots \\
      \frac{1}{n} & \frac{1}{n} & \frac{1}{n}  & \frac{1}{n} & \cdots & \frac{1}{n} \\
    \end{bmatrix}.$$
    $Q$ is also an SIA matrix since the first column of $Q^2$ is positive. However,
     $$QP=\begin{bmatrix}
           \frac{1}{2} & \frac{1}{2} & 0 & 0 & \cdots & 0 \\
       \frac{1}{2} & \frac{1}{2} & 0 & 0 & \cdots & 0 \\
      0 & 0 &1 &  0 & \cdots & 0  \\
+ & + & + & + & \cdots & + \\
      \vdots & \vdots & \vdots & \vdots & \ddots& \vdots \\
      + & + & + & + & \cdots & + \\
    \end{bmatrix},$$
    where ``$+$" denotes an element that is  positive. For two disjoint nonempty sets $\mathcal A=\{1,2\},\tilde{\mathcal A}=\{3\}$, $F_{QP}^k(\mathcal A)=\mathcal A$ and $F_{QP}^k(\tilde{\mathcal A})=\tilde{\mathcal A}$ for any positive integer $k$ and hence $QP$ is not an SIA matrix. $P$ is a stochastic Sarymsakov matrix not in $\mathcal G$. This completes the proof.
 \hfill $\Box$

\subsection{A set  closed under matrix multiplication}
In this subsection, we construct a subset of $\mathcal S$, which is closed under matrix multiplication. This subset consists of the set $\mathcal S_1$ and one specific matrix in $\mathcal V_2$.

Let $R$ be a matrix in $\mathcal V_2$ and satisfies that for any disjoint nonempty sets $\mathcal{A},\ \tilde{\mathcal{A}}\subseteq\mathcal N$, either
\begin{equation}\label{eq:W1}
F_R(\mathcal{A})\cap F_R(\tilde{\mathcal{A}})\neq\emptyset,
\end{equation}
or
\begin{equation}\label{eq:W2}
F_R(\mathcal{A})\cap F_R(\tilde{\mathcal{A}})=\emptyset \text{ and }
|F_R(\mathcal{A})\cup
F_R(\tilde{\mathcal{A}})|\geq|\mathcal{A}\cup\tilde{\mathcal{A}}|.\end{equation}
Such a matrix exists. An example is
\begin{equation}\label{eq:ex2}
R=\begin{bmatrix}
    \frac{1}{n} & \frac{1}{n} & \frac{1}{n} &  \frac{1}{n} & \cdots & \frac{1}{n} \\
    1 & 0 & 0 & 0 &\cdots& 0 \\
    0 & 1 & 0 & 0 &\cdots & 0 \\
    \frac{1}{n} & \frac{1}{n} & \frac{1}{n} & \frac{1}{n} & \cdots &\frac{1}{n}\\
    \vdots & \vdots & \vdots & \vdots & \ddots & \vdots \\
    \frac{1}{n} & \frac{1}{n} & \frac{1}{n} &\frac{1}{n} & \cdots & \frac{1}{n}\\
  \end{bmatrix}.
\end{equation}
To verify that $R$ satisfies this condition, we only have to consider the pair of sets $\mathcal{A}=\{2\},\ \tilde{\mathcal{A}}=\{3\}$ since for other pairs of $\mathcal{A},\ \tilde{\mathcal{A}}$, $F_R(\mathcal{A})\cap F_R(\tilde{\mathcal{A}})\neq\emptyset$. It is clear that $|F_R(2)\cup
F_R(3)|=|\{1,2\}|=|\mathcal{A}\cup\tilde{\mathcal{A}}|$ and $F^2_R(2)\cap
F^2_R(3)=\{1\}$.

\begin{theorem}\label{thm:Sprime}
Suppose that $R$ is a matrix in $\mathcal V_2$ and satisfies that for any disjoint nonempty sets $\mathcal{A},\ \tilde{\mathcal{A}}\subseteq\mathcal N$, (\ref{eq:W1}) or (\ref{eq:W2}) holds. Then, the set $\mathcal S_1^\prime=\mathcal S_1\cup \{R\}$ is closed under matrix multiplication and a compact subset of $\mathcal S_1^\prime$ is a consensus set.
\end{theorem}

{\it Proof:} Let $P$ be a matrix in $\mathcal S_1$. We first show that $RP,PR\in\mathcal S_1$. Given two disjoint nonempty sets $\mathcal{A},\ \tilde{\mathcal{A}}\subseteq\mathcal N$, assume that $F_{RP}(\mathcal{A})\cap F_{RP}(\tilde{\mathcal{A}})=\emptyset$. Since $F_{RP}(\mathcal{A})=F_P(F_R(\mathcal A))$ and $F_{RP}(\mathcal{\tilde{A}})=F_P(F_R(\mathcal {\tilde{A}}))$ based on Lemma \ref{lm:consequent}, one has that $F_R(\mathcal{A})\cap F_R(\tilde{\mathcal{A}})=\emptyset$. In view of the fact that $P$ is a Sarymsakov matrix, one has
\begin{eqnarray*}
|F_{RP}(\mathcal{A})\cup F_{RP}(\tilde{\mathcal{A}})| &=& |F_{P}(F_R(\mathcal{A}))\cup F_{P}(F_{R}(\tilde{\mathcal{A}}))|\\
&>&|F_R(\mathcal{A})\cup
F_R(\tilde{\mathcal{A}})| \\
&\geq& |\mathcal{A}\cup\tilde{\mathcal{A}}|.
\end{eqnarray*}
It follows that $RP$ is a Sarymsakov matrix. Consider the matrix $PR$. Suppose that $\mathcal{A},\ \tilde{\mathcal{A}}\subseteq\mathcal N$ are two disjoint nonempty sets satisfying that $F_{PR}(\mathcal{A})\cap F_{PR}(\tilde{\mathcal{A}})=\emptyset$. One similarly derives that
\begin{eqnarray*}
|F_{PR}(\mathcal{A})\cup F_{PR}(\tilde{\mathcal{A}})|&=&|F_{R}(F_P(\mathcal{A}))\cup F_{R}(F_{P}(\tilde{\mathcal{A}}))|\\
&\geq&|F_P(\mathcal{A})\cup
F_P(\tilde{\mathcal{A}})|\\
&>&|\mathcal{A}\cup\tilde{\mathcal{A}}|.
\end{eqnarray*}
Therefore, $PR$ is a Sarymsakov matrix.

We next show that $R^2\in\mathcal S_1$. Since $R\in\mathcal V_2$, for any disjoint nonempty sets $\mathcal{A},\ \tilde{\mathcal{A}}\subseteq\mathcal N$, there exists an integer $k\leq2$ such that either
\begin{equation}\label{eq:Sprime1}
F^k_R(\mathcal{A})\cap F^k_R(\tilde{\mathcal{A}})\neq\emptyset
\end{equation}or
\begin{equation}\label{eq:Sprime2}
F^k_R(\mathcal{A})\cap F^k_R(\tilde{\mathcal{A}})=\emptyset,\ \mbox{and }|F^k_R(\mathcal{A})\cup
F^k_R(\tilde{\mathcal{A}})|>|\mathcal{A}\cup\tilde{\mathcal{A}}|.
\end{equation}
When (\ref{eq:Sprime1}) holds, it follows from Lemma \ref{lm:consequent} that $F_{R^2}(\mathcal{A})\cap F_{R^2}(\tilde{\mathcal{A}})\neq\emptyset$. When (\ref{eq:Sprime2}) holds, suppose that $F_{R^2}(\mathcal{A})\cap F_{R^2}(\tilde{\mathcal{A}})=\emptyset$. If (\ref{eq:Sprime2}) holds for $k=1$, then from the assumption on $R$, we have
 $$|F_{R^2}(\mathcal{A})\cup F_{R^2}(\tilde{\mathcal{A}})|\geq|F_R(\mathcal{A})\cup
F_R(\tilde{\mathcal{A}})|>|\mathcal{A}\cup\tilde{\mathcal{A}}|;$$
if (\ref{eq:Sprime2}) holds for $k=2$, then we immediately have that $|F_{R^2}(\mathcal{A})\cup F_{R^2}(\tilde{\mathcal{A}})|>|\mathcal{A}\cup\tilde{\mathcal{A}}|.$ Hence, $R^2\in\mathcal S_1$.

Recall that the product of matrices in $\mathcal S_1$ still lies in $\mathcal S_1$. It is clear that $P_2P_1$ is a Sarymsakov matrix for $P_1,P_2\in\mathcal S_1^\prime$. By induction $P_k\cdots P_2P_1\in \mathcal S_1$ for $P_i\in\mathcal S_1^\prime,\ i=1,\ldots,k,$ and any integer $k\geq2$, implying that $\mathcal S_1^\prime$ is closed under matrix multiplication. Then, it immediately follows from Theorem \ref{thm:products1} (5) that a compact subset of $\mathcal S_1^\prime$ is a consensus set. \hfill $\Box$

For a set consisting of the Sarymsakov class and two or more matrices which belong to $\mathcal V_2$ and satisfy that for any disjoint nonempty sets $\mathcal{A},\ \tilde{\mathcal{A}}\subseteq\mathcal N$, (\ref{eq:W1}) or (\ref{eq:W2}) holds, whether it is closed under matrix multiplication depends on those specific matrices in  $\mathcal V_2$.

\subsection{Pattern-symmetric matrices}
In this subsection, we discuss the SIA index of a class of $n\times n$ stochastic matrices, each element $P$ of which satisfies the following pattern-symmetric condition
\begin{equation}\label{eq:typesymmetry}
p_{ij}>0\Leftrightarrow p_{ji}>0, \mbox{for } i\neq j.
 \end{equation}
System (\ref{eq:system}) with bidirectional interactions between agents induces a system matrix satisfying (\ref{eq:typesymmetry}), which arises often in the literature.

We present the following lemma regarding the property of a matrix satisfying (\ref{eq:typesymmetry}).

\begin{proposition}\label{pro:typesymmetry}
Suppose that $P$ satisfies (\ref{eq:typesymmetry}) and is an SIA matrix. Then,
\begin{enumerate}
\item $P\in\mathcal S_2$;
\item if $P$ is symmetric, then $P\in\mathcal S_1$.
\end{enumerate}
\end{proposition}

{\it Proof:} (1) Suppose the contrary holds. Then, there exist two disjoint nonempty sets $\mathcal A,\mathcal {\tilde{A}}\subseteq\mathcal N$, such that
$$F_P^2(\mathcal{A})\cap F_P^2(\tilde{\mathcal{A}})=\emptyset\ \mbox{ and } |F_P^2(\mathcal{A})\cup F_P^2(\tilde{\mathcal{A}})|\leq|\mathcal{A}\cup\tilde{\mathcal{A}}|.$$
Since (\ref{eq:typesymmetry}) holds, one can conclude that for any nonempty set $\mathcal C\subseteq \mathcal N$, $\mathcal C\subseteq F_P^2(\mathcal C)$, implying that
$|F_P^2(\mathcal{A})\cup F_P^2(\tilde{\mathcal{A}})|\geq|\mathcal{A}\cup\tilde{\mathcal{A}}|.$
It follows that $|F_P^2(\mathcal{A})\cup F_P^2(\tilde{\mathcal{A}})|=|\mathcal{A}\cup\tilde{\mathcal{A}}|.$ Then, $F_P^2(\mathcal{A})=\mathcal{A}$ and $F_P^2(\tilde{\mathcal{A}})=\tilde{\mathcal{A}}$, which implies that $F_P^k(\mathcal{A})\cap F_P^k(\tilde{\mathcal{A}})=\emptyset$ for any integer $k$. This contradicts the fact that $P$ is an SIA matrix in view of Lemma \ref{lm:SIAorder2}. We then have $P\in\mathcal S_2$.

(2) If $P\not\in\mathcal S_1$, there exist two disjoint nonempty sets $\mathcal A,\mathcal {\tilde{A}}\subseteq\mathcal N$, such that
$$F_P(\mathcal{A})\cap F_P(\tilde{\mathcal{A}})=\emptyset\ \mbox{and } |F_P(\mathcal{A})\cup F_P(\tilde{\mathcal{A}})|\leq|\mathcal{A}\cup\tilde{\mathcal{A}}|.$$
Since for any set  $\mathcal C\subseteq \mathcal N$,
$$\sum_{i\in\mathcal C,j\in F_P(\mathcal{C})}p_{ij}=|\mathcal C|=\sum_{i\in\mathcal C,j\in F_P(\mathcal{C})}p_{ji}\leq|F_P(\mathcal{C})|,$$
one can only have that $|F_P(\mathcal{A})|=|\mathcal{A}|$ and $|F_P(\tilde{\mathcal{A}})|=|\tilde{\mathcal{A}}|$. This implies
$$\sum_{i\in\mathcal A,j\in F_P(\mathcal{A})}p_{ji}=|F_P(\mathcal{A})|.$$ Combined with the fact that $\mathcal A\subseteq F_P^2(\mathcal A)$, one has that $F_P^2(\mathcal{A})=\mathcal A$. Similarly $F_P^2(\mathcal{\tilde{A}})=\mathcal{\tilde{A}}$. One then has that $F_P^k(\mathcal{A})\cap F_P^k(\tilde{\mathcal{A}})=\emptyset$ for any integer $k$, contradicting the fact that $P$ is SIA. Hence $P\in\mathcal S_1$. \hfill $\Box$

For an SIA and nonsymmetric matrix $P$ satisfying (\ref{eq:typesymmetry}), $P$ is not necessarily a Sarymsakov matrix. An example of such a $P$ is
$$P=\begin{bmatrix}
      0 & 1 & 0 & 0 \\
      \frac{1}{2} & 0 &\frac{1}{2}& 0 \\
      0 & \frac{1}{3} & \frac{1}{3}& \frac{1}{3} \\
      0 & 0 & 1 & 0 \\
    \end{bmatrix}.$$
$P$ is not a Sarymsakov matrix, but $P\in\mathcal S_2$.



With the knowledge of the SIA index, the condition for a consensus set of stochastic symmetric matrices in the literature can be derived directly.
It has been established in Example 7 in \cite{BlOl14} that a compact set $\mathcal P$ of stochastic symmetric matrices is a consensus set if and only if $P$ is an SIA matrix for every $P\in\mathcal P$. The necessity part holds for any  consensus set. As we know from  Proposition \ref{pro:typesymmetry},  a stochastic symmetric matrix $P$ is SIA if and only if $P$ is a Sarymsakov matrix.  The sufficient part becomes clear as the Sarymsakov class is closed under matrix multiplication.

The above claim for stochastic symmetric matrices cannot be extended to stochastic matrices that satisfy (\ref{eq:typesymmetry}). The reason is that a  stochastic matrix satisfying (\ref{eq:typesymmetry}), is not necessarily a Sarymsakov matrix. Hence, in view of Theorem \ref{thm:Sarymlargest}, the product of two such matrices may not be SIA anymore. An example to illustrate this is a set $\mathcal P$ consisting of  two matrices
\begin{align*}
P_1=\begin{bmatrix}
      0 & 1 & 0 & 0 \\
      \frac{1}{2} & 0 & \frac{1}{2} & 0 \\
      0 & \frac{1}{3} & \frac{1}{3} & \frac{1}{3} \\
      0 & 0 & 1 & 0 \\
    \end{bmatrix},\ P_2=\begin{bmatrix}
                          0 &\frac{1}{2} & 0 & \frac{1}{2}\\
                          1 & 0 & 0 & 0 \\
                          0 & 0 & 0 & 1 \\
                          \frac{1}{3} & 0 & \frac{1}{3} & \frac{1}{3} \\
                        \end{bmatrix}.
                        \end{align*}
$P_1,P_2$ both satisfy (\ref{eq:typesymmetry}). However, the matrix product $(P_1P_2)^k$ does not converge to a rank-one matrix as $k\rightarrow\infty$.

\section{A class of generalized Sarymsakov matrices}\label{se:W}
We have seen in Theorem \ref{thm:Sprime}  that the Sarymsakov class plus one specific SIA matrix can lead to a closed set under matrix multiplication that contains $\mathcal S_1$. The property (\ref{eq:W2}) of the matrix $R$ turns out to be critical and we next consider a class of matrices containing all such matrices.

\begin{definition}
A stochastic matrix $P$ is said to belong to a set $\mathcal W$ if for any two disjoint nonempty sets $\mathcal{A},\ \tilde{\mathcal{A}}\subseteq\mathcal N$, either
(\ref{eq:W1}) or (\ref{eq:W2}) with $R$ replaced by $P$ holds.
\end{definition}

The definition of the set $\mathcal W$ relaxes that of the Sarymsakov class by allowing the inequality in (\ref{eq:Sarym2})
not to be strict. It is obvious that $\mathcal S_1$ is a subset of $\mathcal W$.  In addition, $\mathcal W$ is a set that is closed under matrix multiplication. To see this,  we show $PQ\in\mathcal W$ for $P,Q\in\mathcal W$. For any two disjoint nonempty sets $\mathcal{A},\ \tilde{\mathcal{A}}\subseteq\mathcal N$, suppose that
$F_{PQ}(\mathcal{A})\cap F_{PQ}(\tilde{\mathcal{A}})=\emptyset$. It follows from (\ref{eq:W2}) that
\begin{eqnarray*}
|F_{PQ}(\mathcal{A})\cup F_{PQ}(\tilde{\mathcal{A}})|
&=&|F_Q(F_P(\mathcal{A}))\cup F_Q(F_P(\tilde{\mathcal{A}}))|\\
&\geq&|F_P(\mathcal{A})\cup F_P(\tilde{\mathcal{A}})|\\
&\geq& |\mathcal{A}\cup\tilde{\mathcal{A}}|,
\end{eqnarray*}
which implies that $PQ\in\mathcal W$.

Compared with the definition of $\mathcal S_1$, the subtle difference in the inequality in (\ref{eq:W2}) drastically changes the property of $\mathcal W$.  A matrix in $\mathcal W$ is not necessarily SIA. For example, permutation matrices belong to the class $\mathcal W$ since for any disjoint nonempty sets $\mathcal{A},\ \tilde{\mathcal{A}}\subseteq\mathcal N$,
\begin{equation}
F_P(\mathcal{A})\cap F_P(\tilde{\mathcal{A}})=\emptyset \text{ and }
|F_P(\mathcal{A})\cup
F_P(\tilde{\mathcal{A}})|=|\mathcal{A}\cup\tilde{\mathcal{A}}|.\end{equation}
One may expect that the set $\mathcal W\cap\mathcal S$ is closed under matrix multiplication. However, the claim is false and an example to show this is the following two SIA matrices
$$P_1=\begin{bmatrix}
        \frac{1}{3} & \frac{1}{3} & \frac{1}{3} \\
        1 & 0 & 0 \\
        0 & 1 &0\\
      \end{bmatrix},\ P_2=\begin{bmatrix}
      0 & 1 &0 \\
      0 & 0 &1 \\
      \frac{1}{3} & \frac{1}{3} & \frac{1}{3} \\
      \end{bmatrix},$$ where $$P_1P_2=\begin{bmatrix}
        + & + & + \\
       0& 1 & 0 \\
        0 & 0& 1\\
      \end{bmatrix}$$ is  not SIA anymore.

Instead of looking at whether a subset of $\mathcal W$ is a consensus set which concerns the convergence of the matrix product formed by an arbitrary matrix sequence from the subset, we explore the sufficient condition for the convergence of the matrix product of the elements from $\mathcal W$ and its application to doubly stochastic matrices.

\subsection{Sufficient conditions for consensus}

\begin{theorem}\label{thm:W}
Let $\mathcal P$ be a compact subset of $\mathcal W$ and let $P(1),\ P(2),\ P(3),\ldots$ be a sequence of matrices from $\mathcal{P}$. Suppose that $j_1,j_2,\ldots$ is an infinite increasing sequence of the indices such that $P(j_1),P(j_2),\ldots$ are Sarymsakov matrices and $\cup_{r=1}^\infty \{P(j_r)\}$ is a compact set. Then, $P(k)\cdots P(1)$ converges to a rank-one matrix as $k\rightarrow\infty$ if there exists an integer $T$ such that $j_{r+1}-j_r\leq T$ for all $r\geq1$.
\end{theorem}

{\it Proof:} Let $k_0$ be an integer such that $(k_0-1)T+1\geq j_1$. Since $j_{r+1}-j_r\leq T$ for all $r\geq1$, one has that for any integer $k\geq k_0$, the matrix sequence $P(kT),\ldots,P((k-1)T+1)$ contains at least one Sarymsakov matrix, i.e., there exists an integer $i_k$ depending on $k$, $1\leq i_k\leq T$, such that $P((k-1)T+i_k)\in\mathcal S_1$.

We prove that for every integer $k\geq k_0$, $P(kT)\cdots P((k-1)T+1)$ is a Sarymsakov matrix. We only consider those pairs of disjoint nonempty sets $\mathcal{A},\ \tilde{\mathcal{A}}$ satisfying
$$F_{P(kT)\cdots P((k-1)T+1)}(\mathcal{A})\cap
F_{P(kT)\cdots P((k-1)T+1)}(\tilde{\mathcal{A}})=\emptyset.$$
Since $P((k-1)T+i_k)\in\mathcal S_1$, combining with the property of the class $\mathcal W$, it follows that
\begin{equation*}
\begin{split}
&|F_{P(kT)\cdots P((k-1)T+1)}(\mathcal{A})\cup
F_{P(kT)\cdots P((k-1)T+1)}(\tilde{\mathcal{A}})|\\
\geq&|F_{P(kT)\cdots P((k-1)T+2)}(\mathcal{A})\cup
F_{P(kT)\cdots P({(k-1)T+2})}(\tilde{\mathcal{A}})|\geq\cdots\\
\geq&|F_{P(kT)\cdots P({(k-1)T+i_k})}(\mathcal{A})\cup
F_{P(kT)\cdots P({(k-1)T+i_k})}(\tilde{\mathcal{A}})|\\
>&|F_{P(kT)\cdots P({(k-1)T+i_k+1})}(\mathcal{A})\cup
F_{P(kT)\cdots P({(k-1)T+i_k+1})}(\tilde{\mathcal{A}})|\\
\geq&\cdots\geq|\mathcal{A}\cup\tilde{\mathcal{A}}|.
\end{split}
\end{equation*} Therefore, $P(kT)\cdots P((k-1)T+1)$ is a Sarymsakov matrix.

Define  $\mathcal P^\ast=\cup_{r=1}^\infty \{P(j_r)\}$ and $\mathcal Q_T=\{P_T\cdots P_2P_1 | P_i\in\mathcal P,\ i=1,\ldots,T,\ P_s\in\mathcal P^\ast\ \mbox{for some }\\1\leq s\leq T\}$. $\mathcal Q_T$ is a compact set since both $\mathcal P$ and $\mathcal P^\ast$ are compact. Note that the matrices $P(kT)\cdots P((k-1)T+1)\in\mathcal Q_T$ for all $k\geq k_0$ and all matrices  in $\mathcal Q_T$ are Sarymsakov matrices from the above discussion. From Theorem \ref{thm:products1}, we know that
$$\lim_{k\rightarrow\infty}P(kT)P({kT-1})\cdots P((k_0-1)T+1)=\bm{1}c^T,$$ for some nonnegative normalized vector $c$. For any integer $s\geq 1$, there exists an integer $k$ such that $kT+1\leq s<(k+1)T$. Let $||\cdot||$ be the infinity norm of a matrix.  We have
\begin{equation*}
\begin{split}
&||P(s)\cdots P(2)P(1)-\bm{1}c^TP({(k_0-1)T})\cdots P(1)||\\
=&||P(s)\cdots P(1)-P(s)\cdots P({kT+1})\bm{1}c^TP({(k_0-1)T})\cdots P(1)||\\
\leq&||P(s)\cdots P({kT+1})||\cdot||P(kT)\cdots P({(k_0-1)T+1})-\bm{1}c^T||\\
&\cdot||P({(k_0-1)T})\cdots P(1)||\\
\leq&||P(kT)\cdots P({(k_0-1)T+1})-\bm{1}c^T||.
\end{split}
\end{equation*}
Thus, the matrix product $P(s)\cdots P(2)P(1)$ converges to a rank-one matrix as $s$ goes to infinity. \hfill $\Box$


\begin{remark}
In Theorem \ref{thm:W}, if $T_r=j_{r+1}-j_r,\ r\geq1$ is not uniformly upper bounded, we may not be able to draw the conclusion.  The reason is that $\cup_{r=1}^\infty \mathcal Q_{T_r}$ ($\mathcal Q_{T_r}$ can be defined similarly to $\mathcal Q_T$) is not necessarily compact so that the conditions in Theorem \ref{thm:products1} do not apply.
\end{remark}

When the set $\mathcal P$ is a finite set, we have the following corollary.

\begin{corollary}
Let $\mathcal P$ be a finite subset of $\mathcal W$ and let $P(1),\ P(2),\ P(3),\ldots$ be a sequence of matrices from $\mathcal{P}$. Suppose that $j_1,j_2,\ldots$ is an infinite increasing sequence of the indices such that $P(j_1),P(j_2),\ldots$ are Sarymsakov matrices. Then, $P(k)\cdots P(1)$ converges to a rank-one matrix as $k\rightarrow\infty$ if there exists an integer $T$ such that $j_{r+1}-j_r\leq T$ for all $r\geq1$.
\end{corollary}


\subsection{Applications to doubly stochastic matrices}
In fact, the set of matrices $\mathcal W$ contains all doubly stochastic matrices. We can establish the following property of doubly stochastic matrices using the Birkhoff--von Neumann theorem \cite{HoJo85}.

\begin{lemma}\label{lm:doubleSarym1}
Let $P$ be a doubly stochastic matrix. For any nonempty set $\mathcal A\subseteq\mathcal N$, $| F_P(\mathcal A)|\geq|\mathcal A|$.
\end{lemma}

{\it Proof:} From the Birkhoff--von Neumann theorem \cite{HoJo85}, $P$ is doubly stochastic if and only if $P$ is a convex combination of permutation matrices, i.e., $P=\sum_{i=1}^{n!}\alpha_i P_i$ where $\sum_{i=1}^{n!}\alpha_i=1$, $a_i\geq0$ for $i=1,\dots,n!$ and $P_i$ are permutation matrices. For any permutation matrix $P_i$, it is obvious that $| F_{P_i}(\mathcal A)|=|\mathcal A|$ for any set $\mathcal A\subseteq\mathcal N$. In view of the Birkhoff--von Neumann theorem, it holds that
$$ F_P(\mathcal A)=\cup_{\alpha_i\neq0} F_{P_i}(\mathcal A).$$
It then immediately follows that $| F_P(\mathcal A)|\geq|\mathcal A|$.  \hfill $\Box$

From the above lemma, it is easy to see that for a doubly stochastic matrix $P$, either (\ref{eq:W1}) or (\ref{eq:W2}) holds. Hence doubly stochastic matrices belong to the set $\mathcal W$. The following lemma reveals when a doubly stochastic matrix is a Sarymsakov matrix.

\begin{lemma}\label{lm:doubleSarym2}
Let $P$ be a doubly stochastic matrix. $P$ is a Sarymsakov matrix if and only if for every nonempty set $\mathcal A\subsetneq\mathcal N$, $| F_P(\mathcal A)|>|\mathcal A|$.
\end{lemma}

{\it Proof:} The sufficiency part is obvious. We prove the necessity. Suppose on the contrary that there exists a nonempty set $\mathcal A\subsetneq\mathcal N$ such that $|F_P(\mathcal A)|\leq|\mathcal A|$. It follows from Lemma \ref{lm:doubleSarym1} that
\begin{equation}\label{eq:lmdoubleSarym2_1}
|F_P(\mathcal A)|=|\mathcal A|=\sum_{i\in\mathcal A,j\in F_P(\mathcal A)}p_{ij}.
\end{equation}
Since $P$ is doubly stochastic,
\begin{equation}\label{eq:lmdoubleSarym2_2}
\sum_{i\in\mathcal{\bar{A}},j\in F_P(\mathcal A)}p_{ij}=|F_P(\mathcal A)|-\sum_{i\in\mathcal A,j\in F_P(\mathcal A)}p_{ij}=0.
\end{equation}
It follows that $F_P(\bar{\mathcal A})\subseteq\overline{F_P(\mathcal A)}$. Lemma \ref{lm:doubleSarym1} implies that
$$|F_P(\bar{\mathcal A})|\geq|\mathcal{\bar{A}}|=n-|\mathcal A|=|\overline{F_P(\mathcal A)}|.$$
One can conclude that $|F_P(\bar{\mathcal A})|=n-|\mathcal A|$ and
$ F_P(\bar{\mathcal A})=\overline{F_P(\mathcal A)}.$
Then,
$$F_P(\mathcal A)\cap F_P(\bar{\mathcal A})=\emptyset,\ \mbox{and } |F_P(\mathcal A)\cup F_P(\bar{\mathcal A})|=n=|\mathcal A\cup\mathcal{\bar{\mathcal A}}|,$$
which contradicts the fact that $P$ is a Sarymsakov matrix.   \hfill $\Box$


Lemma \ref{lm:doubleSarym2} provides a condition to decide whether a doubly stochastic matrix belongs to $\mathcal S_1$ or not. We have the following corollary based on Theorem \ref{thm:W}.

\begin{corollary}
Let $\mathcal P$ be a compact set of doubly stochastic matrices and let $P(1),\ P(2),\ P(3),\ldots$ be a sequence of matrices from $\mathcal{P}$. Suppose that $j_1,j_2,\ldots$ is an infinite increasing sequence of the indices such that $P(j_1),P(j_2),\ldots$ are Sarymsakov matrices and $\cup_{r=1}^\infty \{P(j_r)\}$ is a compact set. Then, $P(k)\cdots P(1)$ converges to a rank-one matrix as $k\rightarrow\infty$ if there exists an integer $T$ such that $j_{r+1}-j_r\leq T$ for all $r\geq1$.
\end{corollary}

For doubly stochastic matrices with positive diagonals, a sharp statement can be stated.

\begin{proposition}\label{pro:doubleSarym3}
Let $P$ be a doubly stochastic matrix with positive diagonals. If $P$ is SIA, then $P\in\mathcal S_1$.
\end{proposition}

{\it Proof:} Since the diagonal elements of  $P$ are all positive, for any nonempty set $\mathcal A\subseteq\mathcal N$, $\mathcal A\subseteq F_P(\mathcal A)$. Assume that $P\notin\mathcal S_1$. It follows from Lemmas \ref{lm:doubleSarym1} and \ref{lm:doubleSarym2} that there exists a set  $\mathcal A\subsetneq\mathcal N$ such that $| F_P(\mathcal A)|=|\mathcal A|$. One has that $F_P(\mathcal A)=\mathcal A$. By the double stochasticity of $P$, $F_P(\bar{\mathcal A})=\bar{\mathcal A}$, implying that $F_P^k(\mathcal A)\cap F_P^k(\bar{\mathcal A})=\emptyset$ for any integer $k\geq1$. This contradicts the fact that $P$ is SIA.  \hfill $\Box$

For doubly stochastic matrices satisfying condition (\ref{eq:typesymmetry}), we have a similar result.
\begin{proposition}\label{pro:doublesymmetric}
Let $P$ be a doubly stochastic matrix satisfying condition (\ref{eq:typesymmetry}). If $P$ is SIA, then $P\in\mathcal S_1$.
\end{proposition}

{\it Proof:} Assume on the contrary $P$ is not a Sarymsakov matrix. In view of Lemma \ref{lm:doubleSarym2}, there exists a set $\mathcal A\subseteq\mathcal N$ such that $|\mathcal A|=|F_P(\mathcal A)|$. From the proof of Lemma \ref{lm:doubleSarym2} one knows that (\ref{eq:lmdoubleSarym2_1}) and (\ref{eq:lmdoubleSarym2_2}) hold and in addition $|\bar{\mathcal A}|=|F_P(\bar{\mathcal A})|=|\overline{F_P(\mathcal A)}|$. From condition (\ref{eq:typesymmetry}), (\ref{eq:lmdoubleSarym2_1}) implies that  $p_{ij}=0$ for $i\in F_P(\mathcal A),j\in\bar{\mathcal A}$. This implies $F_P^2(\mathcal A)\subseteq\mathcal A$. Combined with the fact that $\mathcal A\subseteq F_P^2(\mathcal A)$, it follows that $\mathcal A=F_P^2(\mathcal A)$. Similarly one has that $\bar{\mathcal A}=F_P^2(\bar{\mathcal A}).$ Therefore $F_P^k(\mathcal A)\cap F_P^k(\bar{\mathcal A})=\emptyset$ contracting the assumption that $P$ is an SIA matrix. \hfill $\Box$

\section{Conclusion}\label{se:conclusion}

In this paper, we have discussed products of generalized stochastic Sarymsakov matrices. With the notion of SIA index, we have shown that the set of all SIA matrices with SIA index no larger than $k$ is closed under matrix multiplication only when $k=1$. Sufficient conditions for the convergence of the matrix product of an infinite matrix sequence to a rank-one matrix have been provided with the help of the Sarymsakov matrices. The results obtained underscore the critical role of the stochastic Sarymsakov class in the set of SIA matrices and in constructing a convergent  matrix sequence to consensus. Construction of a larger set than the one constructed in the paper which is closed under matrix multiplication is a subject for future research.


\bibliographystyle{unsrt}
\bibliography{ref_ming}

\end{document}